\documentclass[11pt,a4paper,leqno]{article}
\usepackage{a4wide}
\usepackage{latexsym}
\usepackage{amsmath}
\usepackage{amssymb}
\usepackage{stackrel}  
\usepackage{color}
\usepackage{graphicx}
\usepackage[linesnumbered,ruled,vlined]{algorithm2e}

\usepackage{etoolbox}
\usepackage{hyperref}

\makeatletter
\newcounter{author}
\renewcommand*\author[1]{%
  \stepcounter{author}%
  \ifnum\c@author=1
    \gdef\@author{#1}%
  \else
    \xdef\@author{\unexpanded\expandafter{\@author\and#1}}%
  \fi
  \csgdef{author@\the\c@author}{#1}}
\newcommand*\email[1]{%
  \csgdef{email@\the\c@author}{#1}}
\newcommand*\address[1]{%
  \csgdef{address@\the\c@author}{#1}}
\AtEndDocument{%
  \xdef\author@count{\the\c@author}%
  \c@author=1
  \print@authors}
\newcommand*\print@authors{%
  \ifnum\c@author>\author@count
  \else
    \print@author{\the\c@author}%
    \advance\c@author by 1
    \expandafter\print@authors
  \fi}
\newcommand*\print@author[1]{%
  \par\medskip
  \begin{tabular}{@{}l@{}}%
    \textsc{\textbf{\csuse{author@#1}}}\\
    \csuse{address@#1}\\
    \textit{E-mail}:
    \href{mailto:\csuse{email@#1}}{\csuse{email@#1}}
  \end{tabular}}
\makeatother

\newtheorem{defin}{Definition}

\newtheorem{theo}{Theorem}
\newtheorem{corol}{Corollary}
\newtheorem{example}{Example}
\newtheorem{remark}{Remark}
\newenvironment{proof}{\medskip\par\noindent{\bf Proof}}{\hfill $\Box$
\medskip\par}

\newcommand{\C}{\mathbb{C}}
\newcommand{\N}{\mathbb{N}}
\newcommand{\Z}{\mathbb{Z}}
\newcommand{\R}{\mathbb{R}}

\newcommand{\Oo}{\mathcal{O}}

\newcommand{\ord}{\mathrm{ord}}
\newcommand{\conv}{\mathrm{conv}}
\newcommand{\coker}{\mathrm{coker}}
\newcommand{\im}{\mathrm{im}}

\newcommand{\corner}{\mathbin{\vrule height 0.13ex depth 0pt width 1.3ex
\vrule height 1.6ex depth 0pt width
0.13ex}\,}

\title{On conditions determining formal automorphisms of integro-differential operators}

\author{Alberto Lastra}
\address{Universidad de Alcal\'a, Departamento de F\'isica y Matem\'aticas, E-28871.\\ 
Alcal\'a de Henares, Madrid, Spain}
\email{alberto.lastra@uah.es}

\author{S{\l}awomir Michalik}
\address{Faculty of Mathematics and Natural Sciences, College of Science,\\Cardinal Stefan Wyszynski Uniwersity in Warsaw, W\'oycickiego 1/3, 01-938 Warszawa, Poland}
\email{s.michalik@uksw.edu.pl}

\author{Maria Suwi{\'n}ska}
\address{Faculty of Mathematics and Natural Sciences, College of Science,\\Cardinal Stefan Wyszynski Uniwersity in Warsaw, W\'oycickiego 1/3, 01-938 Warszawa, Poland}
\email{m.suwinska@op.pl}

\date{}

\begin{document}
\maketitle

\thispagestyle{empty}
{ \small \begin{center}
{\bf Abstract}
\end{center}

Results providing conditions on a family of integro-differential operators to determine a formal automorphism are established. Equivalently, the problem can be read in terms of existence and uniqueness of formal solutions of Cauchy problems in different settings.

The main achievements provide a three-statement result  not only in the framework of general formal power series, but also on subspaces appearing in applications such as Gevrey settings and moment differential operators.

\smallskip

}
\bigskip \bigskip

\section{Introduction}

In the present work, we study conditions for integro-differential operators of the form $P(\partial_t,\partial_z)\partial_{t}^{-m}$ to define a linear automorphism when acting on different spaces of formal power series. Here, $P$ stands for a polynomial in two variables,  $m$ is a positive integer and $\partial_{t}^{-1}$ is the integral operator defined as $\partial_t^{-1}u(t)=\int_0^tu(\tau)\,d\tau$. 

This problem naturally arises when studying the existence and uniqueness of formal solutions to differential Cauchy problems. In this sense, we assume that $P(\partial_t,\partial_z)$ is a general differential operator in two complex variables $(t,z)\in\C^2$ with variable coefficients that belong to the space of formal power series in the variables $(t,z)$, say $\C[[t,z]]$, or in particular to the space of formal power series in the variable $t$, with coefficients being holomorphic functions on some common neighborhood of the origin, say $\Oo[[t]]$. Let us write it in the form 
\begin{equation}
 \label{eq:operatorP}
P(\partial_t,\partial_z)=\sum_{(j,r)\in\Lambda}a_{jr}(t,z)\partial_t^j\partial_z^r,
\end{equation}
where $\Lambda\subset\N_0\times \N_0:=\{0,1,\ldots\}\times\{0,1,\ldots\}$ is a finite set of indices and 
$a_{jr}(t,z)\in\Oo[[t]]$ (resp. $a_{jr}(t,z)\in\C[[t,z]]$) for every $(j,r)\in\Lambda$.

For fixed $m\in\N_0$ one can consider the Cauchy problem
\begin{equation*}
  \left\{
   \begin{array}{l}
    P(\partial_{t},\partial_{z})u(t,z)=f(t,z)\\
     \partial_t^ju(0,z)=\varphi_j(z)\quad\textrm{for}\quad j=0,\dots,m-1
   \end{array}
  \right.
 .\end{equation*}
 
For such Cauchy problem, natural questions arise on the existence and unicity of formal solutions. More precisely, we will be concern on finding conditions on the operator $P(\partial_t,\partial_z)$ given by (\ref{eq:operatorP}), under which statement (A) holds or both statements (A) and (B) hold: 
 \begin{enumerate}
  \item[(A)] For every $f(t,z)\in\C[[t,z]]$ and every $\varphi_j(z)\in\C[[z]]$ ($j=0,\dots,m-1$) there exists exactly one solution $u(t,z)\in\C[[t,z]]$
  \item[(B)] Fix $s\geq 0$. For every $f(t,z)\in \C[[z]]_s[[t]]$ and every  $\varphi_j(z)\in\C[[z]]_s$ ($j=0,\dots,m-1$) there exists exactly one solution $u(t,z)\in \C[[z]]_s[[t]]$.  
 \end{enumerate}

Here, $\C[[z]]_s$ stands for the space of Gevrey formal power series in $z$ of order $s$, and $\C[[z]]_s[[t]]$ denotes the space of formal power series in $t$ with coefficients in the space of Gevrey formal power series $\C[[z]]_s$. The  space of Gevrey series has an essential importance in the development of the classical theory of summability of formal solutions to functional equations. We refer to~\cite{ba2,LodayRichaud} for  further, more in-depth reading on the topic.

Observe that the statement (A) can be reformulated by saying that the integro-differential operator
 \begin{equation}
 \label{eq:operator_A}
  P(\partial_t,\partial_z)\partial_t^{-m}\colon \C[[t,z]]\to\C[[t,z]]
 \end{equation}
is a linear automorphism.

Analogously, given $s\ge0$, the statement (B) can be reformulated by saying that the integro-differential operator
 \begin{equation}
 \label{eq:operator_B}
  P(\partial_t,\partial_z)\partial_t^{-m}\colon \C[[z]]_s[[t]]\to \C[[z]]_s[[t]]
 \end{equation}
is a linear automorphism. 

For this reason, the main results in the present work are focused on giving conditions on $P(\partial_t,\partial_z)$, under which operator (\ref{eq:operator_A}) is an automorphism, and both operators (\ref{eq:operator_A}) and (\ref{eq:operator_B}) are automorphisms. The main results in the present work give answer to these questions in Theorem~\ref{thm3} and Theorem~\ref{thm4}, respectively. The results mentioned previously guarantee that each one of the operators turns out to be an automorphism under an additional condition. The main results are stated in the form of a three-statement result, two of them implying a third one. In this direction, an algebraic and a geometric condition on the problem together imply the operator under study to be a linear automorphism. On the other hand, the fact that the operator is a linear automorphism, together with one of the previous conditions (either the geometric or the algebraic one) imply the remaining condition.

Moreover, we discuss the extensions of the aforementioned characterization to the moment integro-differential operators
\begin{equation*}
P(\partial_{\mathfrak{m}_1,t},\partial_{\mathfrak{m}_2,z})\partial_{\mathfrak{m}_1,t}^{-m},
\end{equation*}
where $\partial_{\mathfrak{m}_1,t}$ and $\partial_{\mathfrak{m}_2,z}$ are moment differential operators defined for given sequences of positive numbers $\mathfrak{m}_1$ and $\mathfrak{m}_2$ (see Section~\ref{sec31} for the definition of moment differential operators).

The study of existence and uniqueness of formal solutions to a given partial differential equation has a significant role in the knowledge of analytic solutions to the problem. Roughly speaking, the theory of summability departs from an existing formal power series solving the problem and constructs the analytic solution by means of a summability procedure, known as Borel-Laplace procedure in its most classical version. 

Therefore, existence of a formal solution to the problem is already significant  by itself. 

On the other hand, the relationship between the analytic and the formal solutions usually arises in the form of an asymptotic expansion. The formal solution is a formal power series with possibly null radius of convergence, and its truncation approximates the analytic solution when working on adequate domains on the complex plane. Typically, these are bounded sectors with vertex at the origin (or the point where the formal asymptotic expansion has been chosen to be performed). 

Uniqueness of the formal solution to the problem is also important as it is the input in the procedure mentioned above.

The strategy is to reduce the problem to linear automorphisms of formal differential equations of one variable of the form
$$P(\partial_z)=\sum_{j\in\Gamma}a_j(z)\partial_z^j,$$
for some polynomial $P$ with holomorphic coefficients near the origin. The result obtained (Theorem~\ref{thm1}) describes a characterization of $P$ to be a linear automorphism on $\C[[z]]$ in terms of geometric properties of the Newton polygon associated with $P$, together with a non-resonance condition on the characteristic polynomial associated with $P$. Another characterization is also obtained when studying linear automorphisms on Gevrey formal power series in one variable (Theorem~\ref{thm2}) in terms of additional geometric properties satisfied by the Newton polygon of $P$. The convergent case follows as a direct corollary. The results are then applied in the two variable setting in Section~\ref{sec4} to achieve the main results (Theorem~\ref{thm3} and Theorem~\ref{thm4}). We conclude the work by linking our results to previously known advances in a different context, and remarking the natural appearance of solutions in the form of formal power series with coefficients being holomorphic functions on a varying neighborhood of the origin. This phenomenon has already been observed in recent researchs with incresing interest among the scientific community, such as \cite{carrillolastra} in the formal settings, and also \cite{tahara3} and \cite{lama4} regarding summability results. In a subsequent research we are planning to study the relation between the sequences of radius of the coefficients of a formal solution, and the operator $P(\partial_t,\partial_z)\partial_t^{-m}$.

It is worth emphasizing that recently there has been a growing interest in such spaces of formal power series with coefficients being holomorphic on different complex neighborhood of the origin, see~\cite{carrillolastra,tahara3,lama4}.

In the last few years, new advances in this theory were achieved in recent studies on the formal solutions to partial differential equations in the complex domain and their summability, such as~\cite{ma22,remy22,yo15}, among many others.

Additionally, knowing upper estimates for the coefficients of the formal solution of a differential problem is useful at the time of determining a Gevrey order of the solution in order to apply an appropriate Borel transform to the formal power series. Therefore, the existence of an automorphism in Gevrey settings provides the knowledge of Gevrey upper bounds associated with the solution, when departing from a known upper bound of the Gevrey order of the forcing term $f$. Such results are known as Maillet-type theorems, and remain also an active field of research of complex partial differential equations. See~\cite{LastraTahara,shirai15} among others, and the references therein. 

We are also exploring the generalization of the main results in the present work in the framework of partial moment differential equations, as mentioned above, due to the relevance that these kind of operators is acquiring in the last decade. Recent results on the summability of formal solutions to such equations are~\cite{lamisu3,lamisu4,mi17}, and also on results of Maillet type in~\cite{lamisu2,su}.

The paper is structured as follows. In Section~\ref{secdefres} we state preliminary definitions and results. It is followed by Section~\ref{sec3}, devoted to the study of the problem in the one variable setting, including several examples of the different situations appearing. A single subsection is taken apart (Section~\ref{sec31}) to describe the more general case of moment differential operators in one variable, and another (Section~\ref{sec32}) focusing on automorphisms under Gevrey settings. In the main section of the present work, Section~\ref{sec4}, we state the main results achieved, Theorem~\ref{thm3} and Theorem~\ref{thm4}.

\textbf{Notation:}

$\N$ stands for the set of positive integers, and $\N_0=\N\cup\{0\}$.

Let $r>0$. $D_r$ stands for the open disc centered at the origin and radius $r$, it is to say, $D_r=\{z\in\C: |z|<r\}$.

Given an open set $U\subseteq\C$, we denote $\mathcal{O}(U)$ the set of holomorphic functions defined in $U$. Given a set $A\subseteq\C$, $C(A)$ stands for the set of continuous functions in $A$ to complex values. The set of formal power series with coefficients in a given Banach space $\mathbb{E}$ (in the variable~$t$) is denoted by $\mathbb{E}[[t]]$. The set of formal power series in $t$ with coefficients being holomorphic functions on some common neighborhood of the origin will be denoted by $\Oo[[t]]$. $\C[[t]]$ stands for the vector space of formal power series with coefficients in $\C$, whereas $\C\{z\}\subseteq\C[[t]]$ represents the subspace of convergent formal power series. We also write $\C[z]$ for the set of polynomials with complex coefficients (in the variable $z$). 

\section{Preliminary definitions and results}\label{secdefres}

In this section, we recall the main definitions and known results to be used in the present work. We mainly work with formal power series with complex coefficients of some Gevrey order $s\ge0$, denoted by $\C[[z]]_s$, which consists of all formal power series $\sum_{n=0}^{\infty}\varphi_nz^n\in\C[[z]]$ such that there exist $C,A>0$ with
$$|\varphi_n|\le C A^n\Gamma(1+sn),\quad n\ge0.$$
Here, $\Gamma(\cdot)$ stands for Gamma function.

%


Let $P\in\Oo[t]$ be a polynomial with holomorphic coefficients (or, in the more general case, with coefficients being formal power series) and consider the formal differential operator $P(\partial_z)$ acting on $\C[[z]]$. Here, $\partial_{z}$ stands for the formal differentiation operator.

\begin{defin}[{see \cite[Section~4]{LodayRichaud}}]
The index $\chi$ of the operator $P(\partial_z)\colon\C[[z]]\to\C[[z]]$ is defined as
\begin{equation*}
\chi(P(\partial_z),\C[[z]]):=\dim\ker(P(\partial_z),\C[[z]])-\dim\coker(P(\partial_z),\C[[z]]),
\end{equation*}
where
\begin{equation*}
\ker(P(\partial_z),\C[[z]]):=\{u\in\C[[z]]:\ P(\partial_z)u(z)=0\},
\end{equation*}
\begin{equation*}
\im(P(\partial_z),\C[[z]])=\{f\in\C[[z]]:\ P(\partial_z)u(z)=f(z)\ \textrm{for some}\ u\in\C[[z]]\}
\end{equation*}
and
\begin{equation*}
\coker(P(\partial_z),\C[[z]])=\C[[z]]/\im(P(\partial_z),\C[[z]]).
\end{equation*}
\end{defin}

Analogously, the index $\chi$ of the operator $P(\partial_z)$ can be defined when restricting the previous operator to the space $\C[[z]]_s$ for any fixed $s\ge0$. By maintaining the same notation for such restriction, we observe that $P(\partial_z) \colon\C[[z]]_s\to\C[[z]]_s$ for any fixed $s\geq 0$ (see for example Proposition 1.2.4,~\cite{LodayRichaud}, together with Stirling's formula).

\section{Linear automorphisms of formal differential operators of one variable}\label{sec3}

In this section, we state equivalent conditions for a formal operator $P(\partial_z):\C[[z]]\to\C[[z]]$ (resp. $P(\partial_z):\C[[z]]_s\to\C[[z]]_s$) to be an automorphism, for some given polynomial $P$ with coefficients given by analytic functions near the origin (resp. some given polynomial $P$ with coefficients given by analytic functions near the origin, and some fixed $s\ge0$). The main results of the present work, dealing with formal differential operators in two variables, lean on those developed in this section.

Let $\Lambda\subset\N_0$ be a finite set of indices and let
\begin{equation}
\label{eq:ordinary_operator}
P(\partial_z)=\sum_{j\in\Lambda}a_j(z)\partial_z^j
\end{equation}
be a differential operator of order $p\in\N$ with holomorphic coefficients on some neighborhood of the origin, say $D$. The order of the zero of $a_j(z)$ at $z=0$ is denoted by $\alpha_j:=\ord_z(a_j)$. Therefore, we may write 
\begin{equation}
 \label{eq:coefficients}
  a_j(z)=\sum_{k\ge \alpha_j}a_{j,k}z^k\in\Oo(D),
\end{equation}
for every $j\in\Lambda$.

Following~\cite{ramis84}, we define the Newton polygon of the operator $P(\partial_z)$ as the convex hull of the union of sets $\corner (j,\alpha_j-j)$ for $j\in\Lambda$,
\begin{equation}
\label{eq:Newton}
 N(P):=\conv\left\{\bigcup_{j\in\Lambda} \corner(j,\alpha_j-j)\right\},
\end{equation}
where $\corner(a,b):=\{(x,y)\in\R^2\colon x\le a,\ y\ge b\}$ denotes the second quadrant of $\R^2$ translated by the vector $(a,b)$.

The number
\begin{equation*}
m:=\min_{j\in\Lambda}(\alpha_j-j)\in\Z
\end{equation*}
denotes the lower ordinate of the Newton polygon $N(P)$ of the operator $P(\partial_z)$.

\begin{remark}
 {\rm Observe that, if $L_m:=\{(x,y)\in\R^2\colon y=m\}$ then
\begin{equation*}
 L_m \cap N(P) = L_m \cap \partial N(P) \subseteq \{(x,y)\in\R^2\colon x\le p,\ y=m\}.
\end{equation*}
}
\end{remark}

\begin{defin}
Given $P$ as before, the principal part $P_{-m}(\partial_z)$ of the operator $P(\partial_z)$ is defined as
\begin{equation}
 \label{eq:principal_part}
 P_{-m}(\partial_z):=\sum_{j\in\Lambda_m}a_{j,j+m}z^{j+m}\partial_z^j,\quad\textrm{where}\quad
 \Lambda_{m}:=\{j\in\Lambda:\ \alpha_j-j=m\}.
\end{equation}

We also define \emph{the characteristic polynomial of the operator $P(\partial_z)$} as
\begin{equation}
\label{eq:characteristic_polynomial}
W_m(\lambda):=\sum_{j\in\Lambda_m}a_{j,j+m}\lambda\cdots(\lambda-j+1).
\end{equation}
\end{defin}

\begin{theo}
\label{thm1}
 The operator $P(\partial_z)$ is a linear automorphism on $\C[[z]]$ if and only if the following conditions hold:
 \begin{enumerate}
  \item[(a)] The lower ordinate $m$ of the Newton polygon $N(P)$ is equal to zero.
  \item[(b)] (non-resonance condition) The characteristic polynomial $W_0(n)$ of the operator $P(\partial_z)$ is different from zero for every $n\in\N_0$.
 \end{enumerate}
\end{theo}

\begin{proof}
We consider the equation
\begin{equation}
 \label{eq:eq}
 P(\partial_z)u(z)=f(z),\quad\textrm{where}\quad u(z)=\sum_{n=0}^{\infty}u_n z^n\in\C[[z]]\quad\textrm{and}\quad f(z)=\sum_{n=0}^{\infty}f_n z^n\in\C[[z]].
\end{equation}
It is worth providing a previous reasoning to the proof that helps on understanding the procedure of determining the existence and/or uniqueness of a formal power series $u(z)\in\C[[z]]$ to satisfy an equation of the form (\ref{eq:eq}), for any given $f(z)\in\C[[z]]$. 

First, observe that the operator $P(\partial_z)$ is a linear automorphism on $\C[[z]]$ if and only if for every $f\in\C[[z]]$ there exists exactly one $u\in\C[[z]]$ satisfying (\ref{eq:eq}), if and only if $\dim\ker(P(\partial_z),\C[[z]])=0$ and $\dim\coker(P(\partial_z),\C[[z]])=0$. 

For $u(z)=\sum_{n=0}^\infty u_n z^n\in\C[[z]]$ we get
\begin{multline}
\label{eq:P_m}
P_{-m}(\partial_z)u(z)=\sum_{j\in\Lambda_m}a_{j,j+m}z^{j+m}\partial_z^j\left(\sum_{n=0}^\infty u_n z^n\right)\\=\sum_{n=0}^\infty\left(\sum_{j\in\Lambda_m}a_{j,j+m}n\cdots(n-j+1)\right)u_n z^{n+m}
=\sum_{n=0}^\infty W_m(n)u_n z^{n+m}.
\end{multline}
On the other hand, we consider the rest of the operator $P(\partial_z)$, notated by $\tilde{P}(\partial_z)$ and defined by
\begin{equation*}
 \tilde{P}(\partial_z):=P(\partial_z)-P_{-m}(\partial_z).
\end{equation*}
We may write this operator as
\begin{equation*}
  \tilde{P}(\partial_z)=\sum_{j\in\tilde{\Lambda}}\tilde{a}_j(z)\partial_z^j
\end{equation*}
for some finite set of indices $\tilde{\Lambda}\subset\N_0$ and some
$\tilde{a}_j(z)\in\Oo(D)$ for $j\in\tilde{\Lambda}$.

It is straightforward to check that $N(\tilde{P})\subset N(P)$ and $N(\tilde{P})\cap L_m = \emptyset$.
More precisely, if we put $\tilde{\alpha}_j:=\ord_z(\tilde{a}_j)$ we conclude that
\begin{equation*}
\tilde{\alpha}_j-j\geq m+1\quad \textrm{for every}\quad j\in\tilde{\Lambda}.
\end{equation*}
Hence, we may write $\tilde{a}_{j}(z)$ as
\begin{equation*}
 \tilde{a}_{j}(z)=\sum_{k=m+1}^{\infty}\tilde{a}_{j,j+k}z^{j+k}.
\end{equation*}
This entails that for every $u(z)=\sum_{n=0}^{\infty}u_n z^n\in\C[[z]]$ we get
\begin{multline}
\label{eq:P_tilde}
\tilde{P}(\partial_z)u(z)=\sum_{j\in\tilde{\Lambda}}\sum_{k=m+1}^{\infty}\tilde{a}_{j,j+k}z^{j+k}\partial_z^j\left(\sum_{n=0}^\infty u_n z^n\right)\\=\sum_{n=0}^\infty\sum_{k=m+1}^{\infty}
\left(\sum_{j\in\tilde{\Lambda}}\tilde{a}_{j,j+k}n\cdots(n-j+1)\right)u_n z^{n+k}
=\sum_{n=0}^\infty\sum_{k=m+1}^{\infty} \tilde{W}_{m,k}(n)u_n z^{n+k},
\end{multline}
where $\tilde{W}_{m,k}(n):=\sum_{j\in\tilde{\Lambda}}\tilde{a}_{j,j+k}n\cdots(n-j+1)$.

If we plug $u(z)=\sum_{n=0}^{\infty}u_nz^n$ into the equation
\begin{equation*}
 P(\partial_z)u(z)=P_{-m}(\partial_z)u(z)+\tilde{P}(\partial_z)u(z)=f(z)=\sum_{n=0}^{\infty}f_n z^n\in\C[[z]],
\end{equation*}
and compare the coefficients at $z^{n+m}$ (for $n\geq 0$ and $n+m\geq 0$) , by (\ref{eq:P_m}) and (\ref{eq:P_tilde}), we conclude that
\begin{multline}
\label{eq:u_n}
 W_m(n)u_n=f_{n+m}-\sum_{k=m+1}^{m+n}\tilde{W}_{m,k}(n-k+m)u_{n-k+m}
 =f_{n+m}-\sum_{j=1}^n\tilde{W}_{m,m+j}(n-j)u_{n-j}\\
 =:A_{m,n}(u_0,\dots,u_{n-1},f_{n+m}),
\end{multline}
where $A_{m,n}\in L(\R^{n+1},\R)$ is a linear form on $\R^{n+1}$.

($\Leftarrow$) If $m=0$ and $W_0(n)\neq 0$ for every $n\in\N_0$ then by (\ref{eq:u_n}) we get
\begin{equation}
\label{eq:u_n_unique}
 u_n=\frac{1}{W_0(n)}A_{0,n}(u_0,\dots,u_{n-1},f_{n})\quad\textrm{for every}\quad n\geq 0.
\end{equation}
It means that for every $f(z)=\sum_{n=0}^{\infty}f_nz^n\in\C[[z]]$ there exists exactly one $u(z)=\sum_{n=0}^{\infty}u_nz^n\in\C[[z]]$ satisfying (\ref{eq:eq}), therefore $\dim\ker(P(\partial_z), \C[[z]])=0$ and $\dim\coker (P(\partial_z),\C[[z]])=0$. Hence, $P(\partial_z)$ is a linear automorphism on $\C[[z]]$.

($\Rightarrow$)
From \cite[Corollary 4.2.5]{LodayRichaud} it follows that the operator $P(\partial_z)\colon\C[[z]]\to\C[[z]]$ has an index with a value  $\chi(P(\partial_z),\C[[z]])=-m$.
Since $P(\partial_z)$ is a linear automorphism on $\C[[z]]$, we conclude that $m=0$.

Assume that the non-resonance condition does not hold and let $n_0\in\N_0$ be the smallest number such that $W_0(n_0)=0$.
We will show that $P(\partial_z)$ is not a surjection.
By (\ref{eq:u_n}), with $m=0$ we get
\begin{equation}
 \label{eq:w_0}
 W_0(n)u_n=f_n-\sum_{k=1}^n\tilde{W}_{0,k}(n-k)u_{n-k}\quad\text{for}\quad n\in\N_0.
\end{equation}
Let $f(z)=\sum_{n=0}^{\infty}f_nz^n=z^{n_0}\in\C[[z]]$.
Then there is no $u(z)=\sum_{n=0}^{\infty}u_nz^n\in\C[[z]]$ such that
$P(\partial_z)u(z)=f(z)$. Indeed, such $u(z)=\sum_{n=0}^{\infty}u_nz^n$ would have to satisfy (\ref{eq:w_0}) for every $n\in\N_0$.
Since $W_0(n)\neq 0$ and $f_n=0$ for $n<n_0$, by (\ref{eq:w_0}) we see that $u_n=0$ for $n<n_0$. Thus, applying (\ref{eq:w_0}) for $n=n_0$ we get the equation $W_0(n_0)u_{n_0}=f_{n_0}$ with $W_0(n_0)=0$ and $f_{n_0}=1$, which is impossible. It means that $f(z)=z^{n_0}\not\in {\rm im}(P(\partial_z),\C[[z]])$, which contradicts the surjectivity of $P(\partial_z)$.
\end{proof}

\begin{remark}
\label{re:coker}
 {\rm If the lower ordinate $m$ of the Newton polygon $N(P)$ is less or equal to zero and the non-resonance condition $W_m(n)\neq0$ holds for every $n\in\N_0$, then by (\ref{eq:u_n}) $\dim\coker(P(\partial_z),\C[[z]])=0$. Consequently, by \cite[Corollary 4.2.5]{LodayRichaud}
 $\dim\ker(P(\partial_z),\C[[z]])=-m$.}
\end{remark}

The reasoning followed in the proof of the previous result can be applied to the following concrete examples. 

\begin{figure}
	\centering
\includegraphics[width=0.45\textwidth]{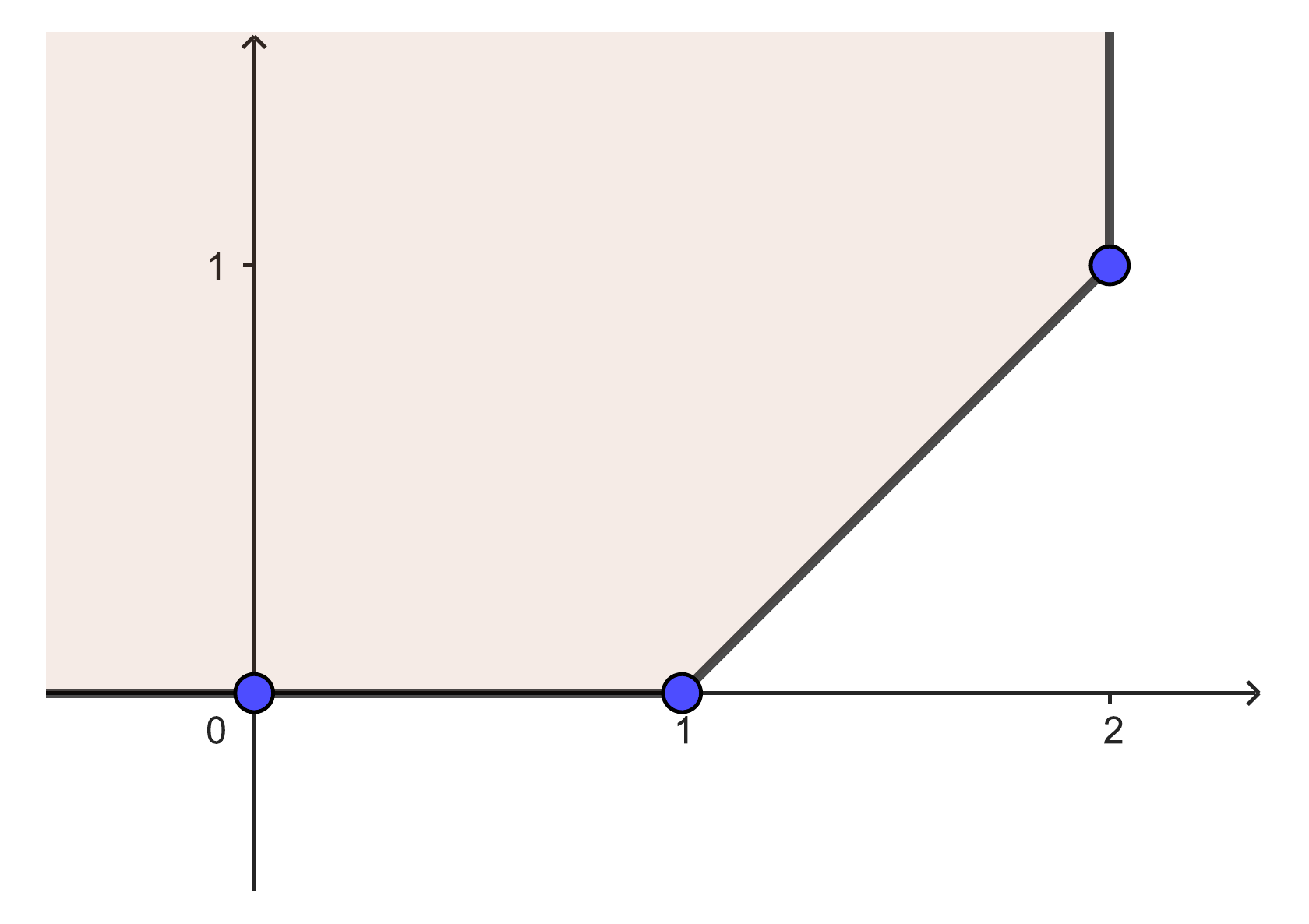}
\includegraphics[width=0.45\textwidth]{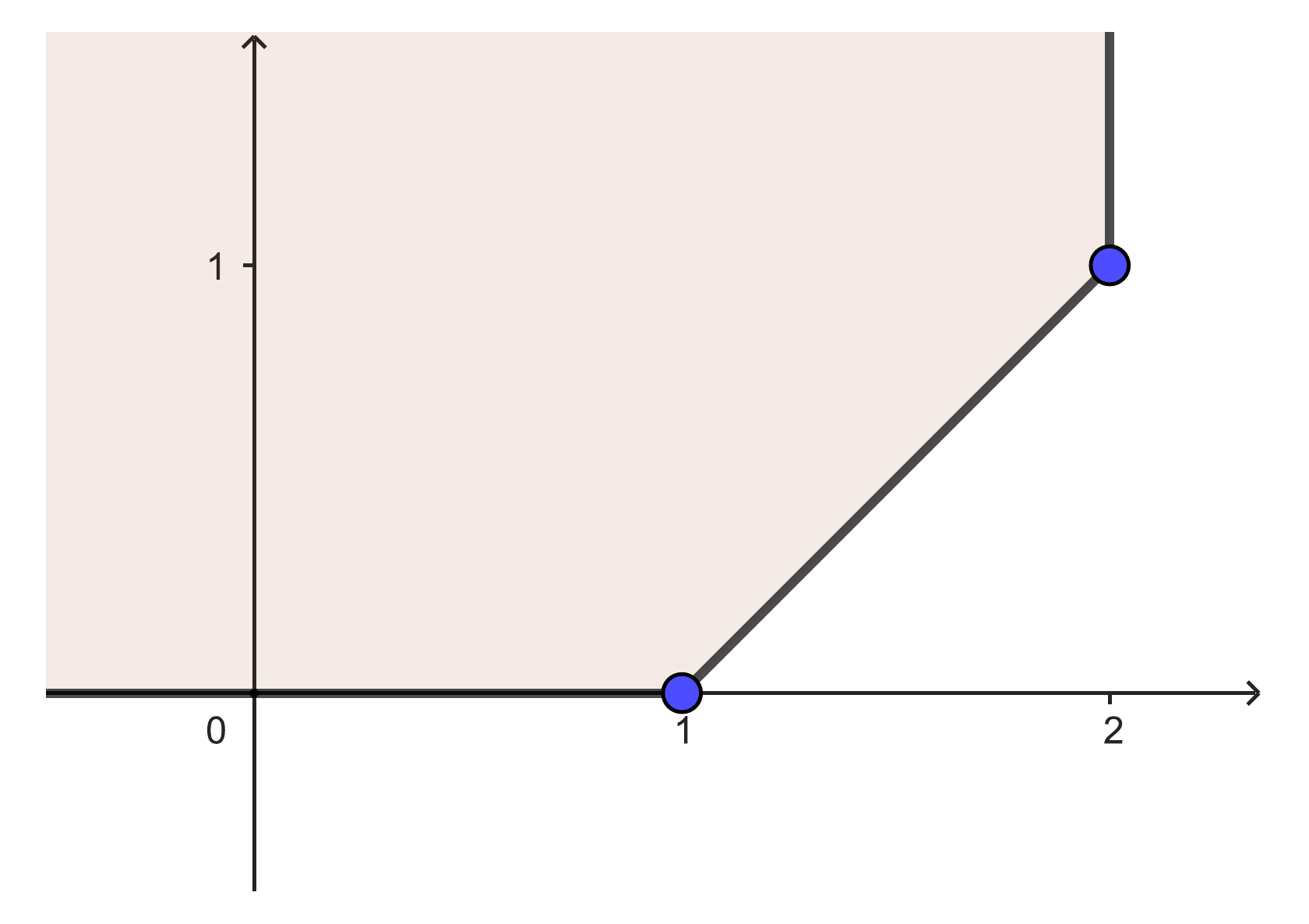}
\caption{Newton polygon associated with $P(\partial_z)$ in Example~\ref{ex1} and~\ref{ex2}}
\label{fig1}
\end{figure}

\begin{example}\label{ex1}
Let $a,b>0$. We consider the operator $P(\partial_z):\C[[z]]\to\C[[z]]$ given by 
$$P(\partial_z)=a+bz\partial_z+z^3\partial_z^2.$$
Its Newton polygon is represented in Figure~\ref{fig1} (left). We observe that $m=0$, so condition (a) in Theorem~\ref{thm1} applies, $\Lambda=\{0,1,2\}$ and $\Lambda_m=\{0,1\}$. It holds that $P_{-m}(\partial_z)=a+bz\partial_z$ and $W_m(\lambda)=a+b\lambda$. Notice that $W_m(n)\neq0$ for any $n\in\N_0$, so condition (b) in Theorem~\ref{thm1} is satisfied. Theorem~\ref{thm1} guarantees that $P(\partial_z)$ is an automorphism of $\C[[z]]$. Given $f(z)=\sum_{n=0}^{\infty}f_nz^n\in\C[[z]]$ the only formal power series $u(z)=\sum_{n=0}^{\infty}u_nz^n$ such that $P(\partial_z)u=f$ is determined as follows: $u_0=f_0/a$, $u_1=f_1/(a+b)$, $u_2=f_2/(a+2b)$, and $u_n=(f_n-(n-1)(n-2)u_{n-1})/(a+bn)$ for $n\ge 3$.
\end{example}

\begin{example}\label{ex2}
Let $b>0$. We consider the operator $P(\partial_z):\C[[z]]\to\C[[z]]$ given by 
$$P(\partial_z)=bz\partial_z+z^3\partial_z^2.$$
Its Newton polygon is represented in Figure~\ref{fig1} (right). We observe that $m=0$, so condition (a) in Theorem~\ref{thm1} applies, $\Lambda=\{1,2\}$ and $\Lambda_m=\{1\}$. It holds that $P_{-m}(\partial_z)=bz\partial_z$ and $W_m(\lambda)=b\lambda$. Notice that $W_m(0)=0$, so condition (b) in Theorem~\ref{thm1} is not satisfied. Theorem~\ref{thm1} guarantees that $P(\partial_z)$ is not an automorphism of $\C[[z]]$. Observe that $P(\partial_z)(c)=0$ for every $c\in\C\subseteq\C[[z]]$, whereas apart from its constant coefficient, given any $f\in\C[[z]]$, all the coefficients of $u\in\C[[z]]$ are determined provided that $P(\partial_z)u=f$. 
\end{example}

\begin{figure}
	\centering
\includegraphics[width=0.45\textwidth]{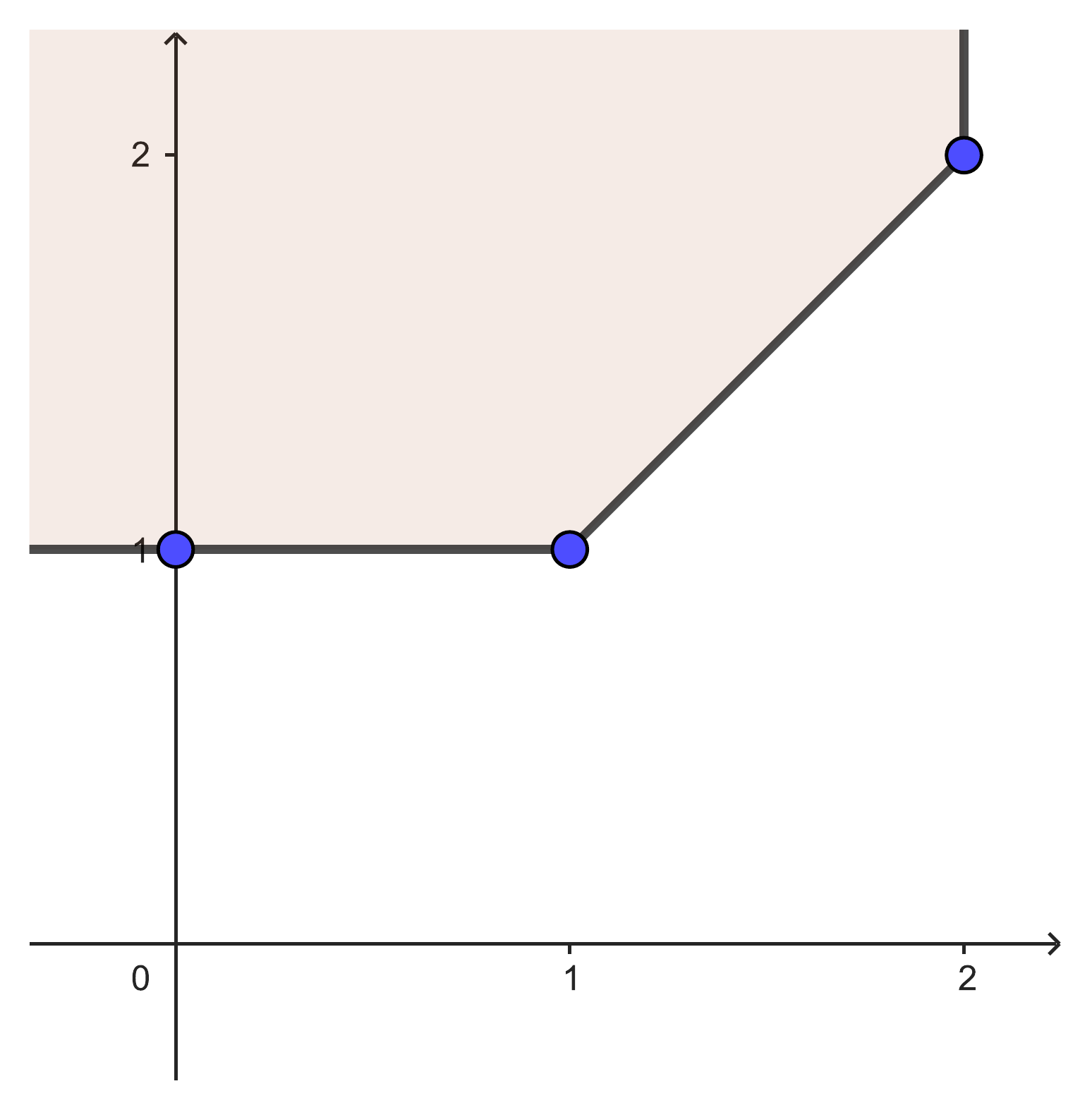}
\includegraphics[width=0.45\textwidth]{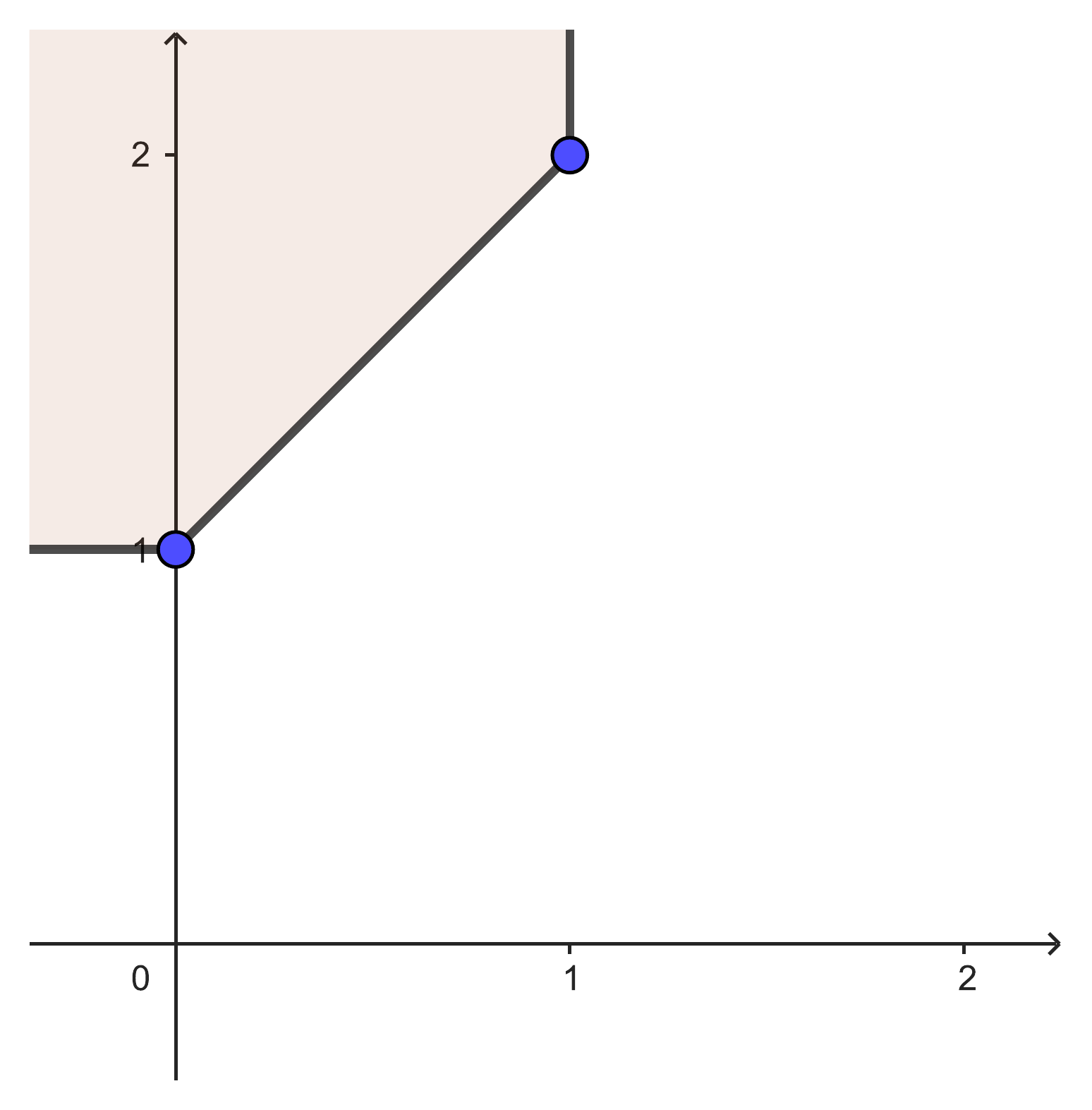}
\caption{Newton polygon associated with $P(\partial_z)$ in Example~\ref{ex3} and~\ref{ex4}}
\label{fig2}
\end{figure}

\begin{example}\label{ex3}
Let $a,b>0$. We consider the operator $P(\partial_z):\C[[z]]\to\C[[z]]$ given by 
$$P(\partial_z)=az+bz^2\partial_z+z^5\partial_z^2.$$
Its Newton polygon is represented in Figure~\ref{fig2} (left). We observe that $m=1$, so condition (a) in Theorem~\ref{thm1} does not hold, $\Lambda=\{0,1,2\}$ and $\Lambda_m=\{0,1\}$. It holds that $P_{-m}(\partial_z)=az+bz^2\partial_z$ and $W_m(\lambda)=a+b\lambda$. Notice that $W_m(n)\neq0$ for all $n\in\N_0$, so condition (b) in Theorem~\ref{thm1} holds. Theorem~\ref{thm1} guarantees that $P(\partial_z)$ is not an automorphism of $\C[[z]]$. Observe that there is no $u(z)\in\C[[z]]$ such that $P(\partial_z)u(z)=1$ as $P(\partial_z)=zP_1(\partial_z)$ for some operator $P_1(\partial_z)$. By direct inspection one can check that, given any element $f(z)\in z\C[[z]]$, there exists a unique $u(z)\in\C[[z]]$ such that $P(\partial_z)u=f$. 
\end{example}

\begin{example}\label{ex4}
Let $a>0$. We consider the operator $P(\partial_z):\C[[z]]\to\C[[z]]$ given by 
$$P(\partial_z)=az+z^3\partial_z.$$
Its Newton polygon is represented in Figure~\ref{fig2} (right). We observe that $m=1$, so condition (a) in Theorem~\ref{thm1} does not hold, $\Lambda=\{0,1\}$ and $\Lambda_m=\{0\}$. It holds that $P_{-m}(\partial_z)=az$ and $W_m(\lambda)=a\lambda$. Notice that $W_m(0)=0$, so condition (b) in Theorem~\ref{thm1} does not hold neither. Theorem~\ref{thm1} guarantees that $P(\partial_z)$ is not an automorphism of $\C[[z]]$. 
\end{example}

\begin{remark}
{\rm In view of (\ref{eq:principal_part}), the condition (a) in Theorem \ref{thm1} means that the principal part of the operator $P(\partial_z)$ defined by (\ref{eq:ordinary_operator}) and (\ref{eq:coefficients}) is a Fuchsian operator given by
 $$
 P_0(\partial_z)=\sum_{j\in\Lambda_0}a_{j,j}z^{j}\partial_z^j.
 $$
 }
\end{remark}

\begin{remark}
{\rm Observe that Theorem \ref{thm1} remains valid when  replacing $a_j(z)\in\Oo(D)$ in (\ref{eq:coefficients}) by $a_j(z)\in\C[[z]]$ for $j\in\Lambda$. A detailed study on the convergence will be precised in Section~\ref{sec32}, when taking $s=0$.}
\end{remark}

\subsection{Linear automorphisms. The moment differential operator setting}\label{sec31}

After minor adaptations Theorem~\ref{thm1} can be generalized to the more general framework of moment derivatives. The importance of moment derivatives in practical applications motivates separating the result in a single subsection, although little details will be given, except for the points, where the argumentation of the proof of Theorem~\ref{thm1} differs.

The concept of moment differentiation was put forward by W. Balser and M. Yoshino in~\cite{bayo}, as a generalization of the classical derivation operator.

\begin{defin}
Let $\mathfrak{m}=(m_n)_{n\ge0}$ be a sequence of positive real numbers. The moment derivative operator $\partial_{\mathfrak{m},z}:\C[[z]]\to\C[[z]]$ is given by
$$\partial_{\mathfrak{m},z}\left(\sum_{n=0}^{\infty}\frac{a_n}{m_n}z^n\right)=\sum_{n=0}^{\infty}\frac{a_{n+1}}{m_n}z^n.$$
\end{defin}

It is known as the moment differential operator due to the fact that the sequence $\mathfrak{m}$ is usually assumed to be a sequence of moments associated with some measure. For example, the classical formal derivative is obtained when considering $\mathfrak{m}$ to be the sequence $(n!)_{n\ge0}$, and moment derivative is closely related to Caputo fractional differential operator when choosing the sequence $\mathfrak{m}$ as $\Gamma_{1/k}:=(\Gamma(1+\frac{n}{k}))_{n\ge0}$, for some fixed $k>0$. Indeed, if one considers Caputo $1/k$-fractional differential operator ${}^{C}D_z^{1/k}$, one has that 
$$(\partial_{\Gamma_{1/k},z}f)(z^{1/k})={}^{C}D_z^{1/k}(f(z^{1/k})),$$
for every $f\in\C[[z]]$. In all the previous cases, the sequence $\mathfrak{m}$ considered is a sequence of moments. Indeed, one has
$$\Gamma(1+\textstyle\frac{p}{k})=\displaystyle\int_{0}^{\infty}s^pks^{k}e^{-s^{k}}ds,\qquad p\in\N_0.$$
A general setting embracing the previous particular examples is that of sequences of moments associated with pairs of kernel functions for generalized summability, developed by J. Sanz in~\cite{sanz}. The construction of Laplace-like operators via the existence of kernel functions is applied in the solution of moment differential equations, see~\cite{jikalasa,lamasa}. 

There are other sequences of great importance in applications, which are closely related to moment sequences. This is the case of the $q-$factorial sequence $\Gamma_{1;q}:=([n]_q!)_{n\ge0}$, for some fixed $q\in\R_{+}\setminus\{1\}$. This sequence is given by  $$[n]_q!=\left\{\begin{aligned}
	1&\quad\textrm{ for }n=0,\\
	[n]_q\cdot[n-1]_q\cdots[1]_q&\quad\textrm{ for any positive integer }n
\end{aligned}\right.$$
with $[j]_q$ standing for the $q-$number $[j]_q=\sum_{h=0}^{j-1}q^{h}$. The $q-$derivative, defined by
$$D_{q,z}f(z)=\frac{f(qz)-f(z)}{qz-z},$$
for every $f\in\C[[z]]$ coincides with the moment derivative $\partial_{\Gamma_{1;q},z}$. When $q>1$, the sequence $\Gamma_{1;q}$ is close to the sequence of moments $(q^{\frac{n(n-1)}{2}})_{n\ge0}$, associated with the Laplace-like operators of kernel $s\mapsto\sqrt{2\pi\ln(q)}\exp\left(\frac{\ln^2(\sqrt{q}s)}{2\ln(q)}\right)$ and also to a kernel involving Jacobi Theta function $s\mapsto\Theta_{1/q}(s)=\sum_{n\in\Z}\frac{1}{q^{\frac{n(n-1)}{2}}}s^n$, both operators appearing in the theory of summability of formal solutions to $q-$difference equations, see~\cite{zhang99,zhang00}.

In principle, $\partial_{\mathfrak{m},z}$ is only defined on formal power series, and consequently on analytic functions near some point by identification of the function with its Taylor expansion at that point. In~\cite{lamisu,lamisu4}, the domain of $\partial_\mathfrak{m}$ has been extended to analytic functions defined on sectors of the complex plane which represent the sum or multisum of a formal asymptotic expansion at the vertex of the sector.

The formal nature of Theorem~\ref{thm1} allows to adapt the arguments of its proof when considering an operator $P(\partial_\mathfrak{m})$, $\mathfrak{m}=(m_n)_{n\ge0}$ being any sequence of positive real numbers normalized by $m_0=1$. 

In the following results, we write $\partial_\mathfrak{m}$ instead of $\partial_{\mathfrak{m},z}$ for simplicity, and we assume $P(\partial_\mathfrak{m})$ is defined by
$$P(\partial_\mathfrak{m})=\sum_{j\in\Lambda}a_{j}(z)\partial_\mathfrak{m}^{j}.$$
Here, $\Lambda\subseteq\N_0$ is a finite set, $a_{j}$ is an analytic function on some neighborhood of the origin (or a~formal power series in $z$) for all $j\in\Lambda$ of the form (\ref{eq:coefficients}). We write $\partial_{\mathfrak{m}}^{0}=\hbox{Id}$, $\partial_\mathfrak{m}^{1}=\partial_\mathfrak{m}$ and $\partial_{\mathfrak{m}}^{j+1}=\partial_\mathfrak{m}\circ\partial_{\mathfrak{m}}^{j}$ recursively for every $j\ge 1$. Therefore, 
$$P(\partial_\mathfrak{m}):\C[[z]]\to\C[[z]].$$

\begin{defin} 
The Newton polygon associated with $P(\partial_\mathfrak{m})$ is defined as the Newton polygon associated with $P(\partial_z)$.
\end{defin}

The lower ordinate of Newton polygon will still be denoted by $m$, and we preserve the definition of $\Lambda_m$. The generalized characteristic polynomial of $P(\partial_{\mathfrak{m}})$, denoted $W_{m,\mathfrak{m}}$ is defined by
$$W_{m,\mathfrak{m}}(\lambda)=\sum_{j\in\Lambda_m,\ j\le \lambda}a_{j,j+m}\frac{m_{\lambda}}{m_{\lambda-j}}.$$

Observe that $W_{m,(n!)_{n\ge0}}(\lambda)$ coincides with $W_{m}(\lambda)$.

The proof of Theorem~\ref{thm1} can be mimicked to prove the following result.

\begin{corol}
Let $P$ be as above and let $\mathfrak{m}=(m_n)_{n\ge0}$ be a sequence of positive real numbers. The operator $P(\partial_\mathfrak{m})$ is a linear automorphism on $\C[[z]]$ if and only if the following conditions hold:
\begin{itemize}
\item[$(a_m)$] The lower ordinate of the Newton polygon $N(P)$ is equal to zero.
\item[$(b_m)$] (non-resonance condition) The generalized characteristic polynomial $W_{0,\mathfrak{m}}(n)$ of $P(\partial_{\mathfrak{m}})$ is different from zero for all $n\in\N_0$.
\end{itemize}
\end{corol}

In the particular case of $\mathfrak{m}=\Gamma_{1/k}$, the operator under study is 
$$P({}^{C}D_z^{1/k})=\sum_{j\in\Lambda}a_j(z^{1/k}){}^{C}D_z^{1/k}:\C[[z^{1/k}]]\to\C[[z^{1/k}]],$$
dealing with automorphisms of $\C[[z^{1/k}]]$.

\vspace{0.3cm}

Considering $\mathfrak{m}=\Gamma_{1;q}$ for some $q>1$, we deal with $P(D_{q,z}):\C[[z]]\to\C[[z]]$, with $W_{m,\Gamma_{1;q}}(\lambda)=\sum_{j\in\Lambda_m,\ j\le \lambda}a_{j,j+m}\prod_{h=\lambda-j+1}^{\lambda}[h]_q$. Observe the confluence of all the constructions under these settings to those studied in the first part of Section~\ref{sec3}.

\subsection{Linear automorphisms. Gevrey settings}\label{sec32}

In this subsection we state equivalent conditions for an automorphism $P(\partial_z)$ of $\C[[z]]$ to remain an automorphism when restricted to $\C[[z]]_{s}$. In particular, we study the convergent case for $s=0$.

We start with the following corollary from the Gevrey Index Theorem for linear differential operators \cite[Corollary 4.2.5]{LodayRichaud} (see also \cite{Malgrange} and \cite{ramis84}).
\begin{theo}
\label{thm2}
 Let $s\geq 0$. The operator $P(\partial_z)$ is a linear automorphism on $\C[[z]]$ which extends to the automorphism on $\C[[z]]_s$  if and only if the following conditions hold:
 \begin{enumerate}
  \item[(a)] The lower ordinate $m$ of the Newton polygon $N(P)$ is equal to zero.
  \item[(b)] (non-resonance condition) $W_0(n)\ne 0$ for every $n\in\N_0$.
	\item[(c)] The first positive slope of the Newton polygon $N(P)$ is greater or equal to $k=1/s$.
 \end{enumerate}
\end{theo}
\begin{proof}

($\Rightarrow$)
Since $P(\partial_z)\colon \C[[z]]\to\C[[z]]$ is a linear automorphism, by Theorem \ref{thm1} we get the assertions (a) and (b). Moreover, since $\chi(P(\partial_z),\C[[z]]_s)=0$, by \cite[Corollary 4.2.5]{LodayRichaud} we conclude that (c) also holds.

($\Leftarrow$) If $m=0$ and $W_0(n)\ne 0$ for every $n\in\N_0$ then by
Theorem \ref{thm1} the operator $P(\partial_z)\colon\C[[z]]\to\C[[z]]$ is a linear automorphism and, consequently, $\dim\ker(P(\partial_z),\C[[z]])=0$. Since $\C[[z]]_s\subset\C[[z]]$, we conclude that $\ker(P(\partial_z),\C[[z]]_s)\subseteq \ker(P(\partial_z),\C[[z]])$, so also $\dim\ker(P(\partial_z),\C[[z]]_s)=0$. If the first positive slope of the Newton polygon $N(P)$ is greater or equal to $k=1/s$,
then by \cite[Corollary 4.2.5]{LodayRichaud} we see that $\chi(P(\partial_z),\C[[z]]_s)=0$, so also $\dim\coker(P(\partial_z),\C[[z]]_s)=0$. Hence, $P(\partial_z)\colon\C[[z]]_s\to\C[[z]]_s$ is a linear automorphism.
\end{proof}

\begin{example}\label{ex5}
Recall that $P(\partial_z)$ determines an automorphism of $\C[[z]]$ for $P(\partial_z)$ given in Example~\ref{ex1}. Theorem~\ref{thm2} guarantees that $P(\partial_z)$ continues to be an automorphism of $\C[[z]]_s$ for every $s\ge1$. This can be proved directly. Indeed, assume that $P(\partial_z)u=f$ for some $u(z)=\sum_{n\ge0}u_nz^n\in\C[[z]]$ and $f(z)=\sum_{n\ge0}f_nz^n\in\C[[z]]_s$ for some $s\ge0$. This entails that $|f_n|\le C_1A_1^n\Gamma(1+sn)\le CA^nn!^s$ for some $C_1,C,A_1,A>0$  by Stirling's formula. 

First assume that $s\ge 1$. Observe from the recursion formula 
$$u_n=\frac{f_n-(n-1)(n-2)u_{n-1}}{a+bn},\qquad n\ge3$$ 
that one can choose sufficiently large constants $\tilde{C},\tilde{A}>0$ such that $|u_0|=|f_0|/a\le\tilde{C}$, $|u_j|=|f_j|/(a+jb)\le \tilde{C}+\tilde{A}^jj^s$ for $j=1,2$ and with $\tilde{C}\ge C$, $\tilde{A}\ge\max\left\{A,\frac{\max\{A,1\}}{\min\{a,b\}}\right\}$. One can also prove by induction that 
\begin{equation}\label{e504}
|u_n|\le \tilde{C}\tilde{A}^nn!^s,\qquad n\ge0.
\end{equation}
Again, Stirling's formula allows to conclude that $u\in\C[[z]]_s$. The estimate (\ref{e504}) is valid for $n=0,1,2$. Assume it is valid up to some $n-1\ge1$. Then,
\begin{multline*}
|u_n|\le\frac{1}{a+bn}\left[|f_n|+(n-1)(n-2)|u_{n-1}|\right]\le\frac{1}{a+bn}\left[CA^nn!^s+n^2\tilde{C}\tilde{A}^{n-1}(n-1)!^s\right]\\
\le n!^s\frac{CA^n+n\tilde{C}\tilde{A}^{n-1}}{a+bn}\le \tilde{C}\tilde{A}^{n-1} n!^s\frac{n+A}{\min\{a,b\}(n+1)}\le\tilde{C}\tilde{A}^nn!^s.
\end{multline*}
This guarantees that $P(\partial_z):\C[[z]]_s\to\C[[z]]_s$ is an automorphism for $s\ge1$. If $0<s<1$, and we consider $f(z)=\sum_{n\ge0}n!^sz^s\in\C[[z]]_s$, one can prove in an analogous manner the existence of constants $\tilde{C}_1,\tilde{A}_1>0$ such that $|u_n|\ge \tilde{C}_1\tilde{A}_1^nn!$, so $P(\partial_z)$ is not an automorphism. 
\end{example}

Taking $s=0$ in Theorem \ref{thm2} we get 
\begin{corol}
 The operator $P(\partial_z)$ of order $p\in\N$ is a linear automorphism on $\C[[z]]$, which extends to the automorphism on $\C\{z\}$ if and only if the following conditions hold:
 \begin{enumerate}
  \item[(a)] $N(P(\partial_z))=\corner(p,0)$,
  \item[(b)] (the non-resonance condition) $W_0(n)\ne 0$ for every $n\in\N_0$.
 \end{enumerate}
\end{corol}

\begin{example}\label{ex6}
Example~\ref{ex1} describes an automorphism of $\C[[z]]$ which does not extend to an automorphism of $\C\{z\}$ in view of the Newton polygon associated with $P(\partial_z)$. See Figure~\ref{fig1} (left).
\end{example}

Let us consider a slightly modified example.

\begin{figure}
	\centering
\includegraphics[width=0.45\textwidth]{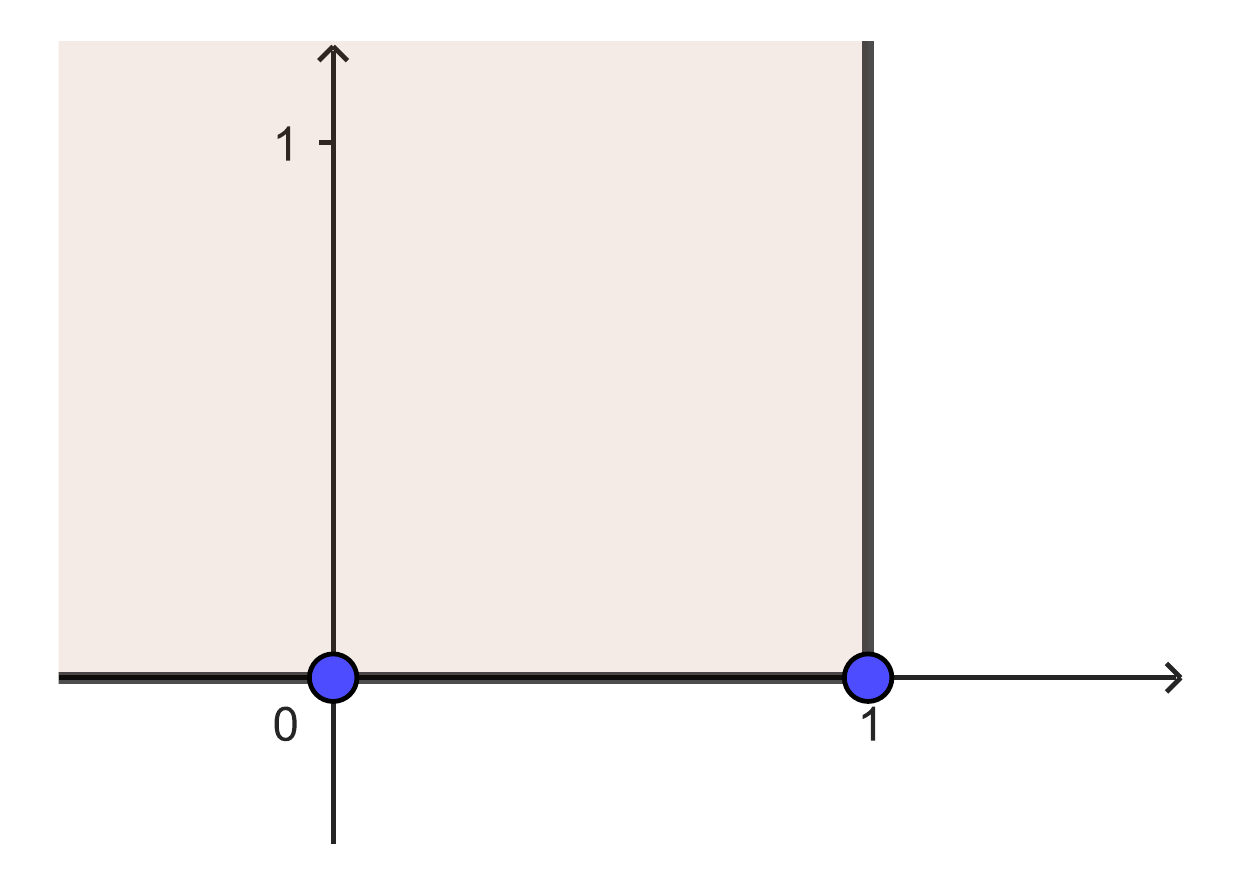}
\caption{Newton polygon associated with $P(\partial_z)$ in Example~\ref{ex5}}
\label{fig3}
\end{figure}

\begin{example}\label{ex7}
Let $a,b>0$. We consider the operator $P(\partial_z):\C[[z]]\to\C[[z]]$ given by
$$P(\partial_z)=a+bz\partial_z.$$
Its Newton polygon is represented in Figure~\ref{fig3}. We observe that $m=0$, so condition (a) in Theorem~\ref{thm1} applies. We observe that $W_m(\lambda)=a+b\lambda$, so that condition (b) in Theorem~\ref{thm1} holds. This entails that $P(\partial_z)$ is an automorphism of $\C[[z]]$. Condition (a) of Corollary~\ref{coro521} also holds. Therefore, $P(\partial_z)$ is also an automorphism of $\C\{z\}$.

Observe that, given $f(z)=\sum_{n\ge0}f_nz^n\in\C[[z]]$, the formal power series $u(z)=\sum_{n\ge0}u_nz^n\in\C[[z]]$ defined by $u_0=f_0/a$, and for all $n\ge1$ $u_n=f_n/(a+bn)$ is the unique formal power series satisfying $P(\partial_z)u=f$. Additionally, if there exist $C,A>0$ such that $|f_n|\le CA^{n}$, for every $n\ge0$, i.e., $f\in\C\{z\}$ with radius of convergence at least $1/A$, then one has that $|u_0|=|f_0|/a\le C/a$ and $|u_n|=|f_n|/(a+bn)\le C/aA^n$ for every $n\ge0$. As a consequence, $u\in\C\{z\}$ with radius of convergence at least $1/A$. Indeed, observe that given $f\in\C\{z\}$ the equation $P(\partial_z)u=au+bz\partial_zu=f(z)$ is a first order linear ODE which can be solved by variation of constants formula for a concrete $f$.
\end{example}

In the last part of this section we showcase how the properties of a differential operator can be derived from a simpler operator maintaining its most significant terms, the principal part of the differential operator. Its detection is now explained.

 By (\ref{eq:ordinary_operator}) and (\ref{eq:coefficients}) we may write the operator $P(\partial_z)$ as
 \begin{equation*}
  P(\partial_z)=\sum_{(j,k)\in\Lambda\times\N_0} a_{jk}z^k\partial_z^j.
 \end{equation*}
Next, for the Newton polygon $N(P)$ given by (\ref{eq:Newton}) we define a set of indices $\Delta\subset \Lambda\times\N_0$ as
\begin{equation*}
 \Delta:=\{(j,k)\in\Lambda\times\N_0\colon a_{jk}\neq 0,\ (j,k-j)\in \partial N(P)\}.
\end{equation*}
We are in position to define the principal part $P_N(\partial_z)$ of the operator $P(\partial_z)$ with respect to the Newton polygon $N(P)$ as
\begin{equation*}
 P_N(\partial_z)=\sum_{(j,k)\in\Delta} a_{jk}z^k\partial_z^j.
\end{equation*}

Observe that the operators $P(\partial_z)$ and $P_N(\partial_z)$ have the same Newton polygon, the same principal part and the same characteristic polynomial. Hence, by Theorems \ref{thm1} and \ref{thm2} we get that the operator $P(\partial_z)$ has the same properties as the operator $P_N(\partial_z)$. More precisely, we have 

\begin{corol}\label{coro521}
 The operator $P(\partial_z)$ is a linear automorphism on $\C[[z]]$ if and only if the operator $P_N(\partial_z)$ is a linear automorphism on the same space.
 
 Moreover, for fixed $s\geq 0$ the automorphism $P(\partial_z)$ on $\C[[z]]$ extends to the automorphism on $\C[[z]]_s$ if and only if the automorphism $P_N(\partial_z)$ extends in the same way.
\end{corol}

The previous result is very important in applications, allowing to simplify certain concrete problems in a significant manner. The following example shows its importance.

\begin{example}
Consider the differential operator 
$$P(\partial_z)=p_1(z)+zp_2(z)\partial_z+z^3\partial_z^2,$$ 
where $p_1,p_2\in\C[z]$ with $p_1(0)=a$ for some $a>0$ and $p_1(0)=b$ for some $b>0$. Then, $P(\partial_z)$ is a linear automorphism on $\C[[z]]$ which extends to a linear automorphism on $\C[[z]]_s$ only for $s\ge1$. This is a consequence of Example~\ref{ex1} and Example~\ref{ex5} and the fact that $P(\partial_z)$ shares its Newton polygon with the differential operator in these examples.
\end{example}

\section{Integro-differential operators of two variables}\label{sec4}
In this section we assume that $\Lambda$ is a finite subset of indices in
$\N_0\times\N_0$. We consider a partial differential operator of the form 
\begin{equation}
\label{eq:operator_p}
P(\partial_t,\partial_z)=\sum_{(j,r)\in\Lambda}a_{jr}(t,z)\partial_t^j\partial_z^r \quad \textrm{with} \quad a_{jr}(t,z)=\sum_{n=\ord_t(a_{jr})}^\infty a_{j,r,n}(z)t^n\in\Oo[[t]]\ (\textrm{or}\ \C[[t,z]]).
\end{equation}
For such operator $P(\partial_t,\partial_z)$ we define the Newton polygon (see \cite{Yonemura}) as the convex hull of the union of sets $\corner (j+r,\ord_t(a_{jr})-j)$ for $(j,r)\in\Lambda$, namely 
\begin{equation}
\label{eq:Newton_2}
 N(P):=\conv\left\{\bigcup_{(j,r)\in\Lambda}\corner(j+r,\,\ord_t(a_{jr})-j)\right\}.
\end{equation}
Let us recall that
$\corner(a,b)=\{(x,y)\in\R^2\colon x\le a,\ y\ge b\}$ denotes the second quadrant of $\R^2$ translated by the vector $(a,b)$.

In this section we assume that 
\begin{equation*}
m:=\max_{(j,r)\in\Lambda}(j-\ord_t(a_{jr}))\ge 0.
\end{equation*}
Observe that the initial problem
\begin{equation*}
  \left\{
   \begin{aligned}
    P(\partial_{t},\partial_{z})u(t,z)&=f(t,z)&\\
     \partial_t^ju(0,z)&=\varphi_j(z)&\ \textrm{for}\quad j=0,\dots,m-1
   \end{aligned}
  \right.
 \end{equation*}
has exactly one solution $u(t,z)\in\Oo[[t]]$ for every $f(t,z)\in\Oo[[t]]$ and $\varphi_j(z)\in\Oo(D)$ ($j=0,\dots,m-1$) if and only if the integro-differential operator
\begin{equation}
\label{eq:integro_differential}
 P(\partial_t,\partial_z)\partial_t^{-m}\colon \Oo[[t]] \to \Oo[[t]]
\end{equation}
is a linear automorphism, where $\partial_{t}^{-1}$ is the integral operator defined as $\partial_t^{-1}u(t)=\int_0^tu(\tau)\,d\tau$. 

For this reason, in this section we reformulate the problem to studying necessary and sufficient conditions, under which the integro-differential operator (\ref{eq:integro_differential}) is a linear automorphism. 
To this end we consider the principal part of the operator $P(\partial_t,\partial_z)$ with respect to $\partial_t$ given by (\ref{eq:operator_p}), which is defined by
\begin{equation}
\label{eq:operator_p_m}
P_m(\partial_t,\partial_z):=\sum_{(j,r)\in\Lambda_m}a_{j,r,j-m}(z)t^{j-m}\partial_t^j\partial_z^r,
\end{equation}
where $\Lambda_m:=\{(j,r)\in\Lambda\colon j-\ord_t(a_{jr})=m\}$.

Before stating the main results of the present work, we provide some preliminary auxiliary constructions.

Since $j\geq m$ for $(j,r)\in\Lambda_m$, one may apply the operator (\ref{eq:operator_p}) to the formal power series
$u(t,z)=\sum_{n=0}^{\infty}u_n(z)t^n$ in $\Oo[[t]]$ (or $\C[[t,z]]$) 
to obtain that
\begin{multline*}
P_m(\partial_t,\partial_z)\partial_t^{-m}\left(\sum_{n=0}^\infty u_n(z)t^n\right)=\\=\sum_{n=0}^\infty\left(\sum_{(j,r)\in\Lambda_m}a_{j,r,j-m}(z)n(n-1)\cdots(n-(j-m)+1)\partial_z^r\right) u_n(z)t^n
=\sum_{n=0}^\infty \tilde{P}_m(n,\partial_z)u_n(z)t^n,
\end{multline*}
where $a_{j,r,j-m}(z)=z^{\alpha_{j,r}}\tilde{a}_{j,r,j-m}(z)$ with $\alpha_{j,r}=\ord_z(a_{j,r,j-m}(z))$, and 
\begin{equation}
\label{eq:operator_p_m_tilde}
\tilde{P}_m(n,\partial_z)=\sum_{(j,r)\in\Lambda_m}\tilde{a}_{j,r,j-m}(z)z^{\alpha_{j,r}}n(n-1)\cdots(n-(j-m)+1)\partial_z^r.
\end{equation}

To construct the principal part of the above operator $\tilde{P}_m(n,\partial_z)$ with respect to $\partial_z$, we put
\begin{equation}
\label{eq:l}
l:=\min_{(j,r)\in\Lambda_m}(\alpha_{j,r}-r)\quad\textrm{and}\quad
\Lambda_{m,l}:=\{(j,r)\in\Lambda_m\colon \alpha_{j,r}-r=l\}.
\end{equation}

Then, for given $n\in\N_0$ the principal part of $\tilde{P}_m(n,\partial_z)$ is given by
\begin{equation*}
 \tilde{P}_{m,-l}(n,\partial_z)=\sum_{(j,r)\in\Lambda_{m,l}}\tilde{a}_{j,r,j-m}(z)z^{r+l}n(n-1)\cdots(n-(j-m)+1)\partial_z^r,
\end{equation*}
under condition that $\tilde{P}_{m,-l}(n,\partial_z)\neq 0$.

By plugging the power series $u_n(z)=\sum_{k=0}^\infty u_{nk}z^k$ in $\C[[z]]$ (or $\C[[z]]_s$ for some $s\geq 0$) into the previous operator, we arrive at
\begin{equation*}
 \tilde{P}_{m,-l}(n,\partial_z)u_n(z)=\sum_{k=0}^\infty\sum_{(j,r)\in\Lambda_{m,l}}\tilde{a}_{j,r,j-m}(z)n(n-1)\cdots(n-(j-m)+1)k(k-1)\cdots(k-r+1)u_{nk}z^{k+l}.
\end{equation*}

Therefore, 
\begin{equation*}
 \tilde{P}_{m,-l}(n,\partial_z)u_n(z)=\sum_{k=0}^\infty W_{m,l}(n,k,z)u_{nk}z^{k+l},
\end{equation*}
where
\begin{equation}
 \label{eq:W_ml}
 W_{m,l}(n,k,z)=\sum_{(j,r)\in\Lambda_{m,l}}\tilde{a}_{j,r,j-m}(z)n(n-1)\cdots(n-(j-m)+1)k(k-1)\cdots(k-r+1).
\end{equation}

At this point, we are in a position to prove the first main result of the present work.

\begin{theo}
\label{thm3}
Any two of the following three conditions entail a third one:
  \begin{enumerate}
 \item[(i)] The operator $P(\partial_t,\partial_z)\partial_t^{-m}$ is a linear automorphism on $\C[[t,z]]$.
  \item[(ii)] The lower ordinate of the Newton polygon $N(\tilde{P}_m(n,\partial_z))$ is equal to $l$ for every $n\in\N_0$, where $l$ is defined by (\ref{eq:l}) and is equal to zero.
  \item[(iii)] The non-resonance condition $W_{m,l}(n,k,0)\ne 0$ holds for every $n,k\in\N_0$.
 \end{enumerate}
\end{theo}

\begin{proof}
The proof is divided into two steps. In a first step, we reduce the problem to that in the one variable settings, already solved in Theorem~\ref{thm1}. In a second step, we adapt the necessary and sufficient conditions obtained in that result to the several variable framework.

First, let us define $Q(\partial_t,\partial_z):=P(\partial_t,\partial_z)-P_m(\partial_t,\partial_z)$, where $P_m(\partial_t,\partial_z)$, given by (\ref{eq:operator_p_m}), is the principal part of the operator $P(\partial_t,\partial_z)$. Then the equation
 \begin{equation*}
  P(\partial_t,\partial_z)\partial_t^{-m}u=f
 \end{equation*}
 for $u=\sum_{n=0}^\infty u_n(z)t^n$ and $f=\sum_{n=0}^\infty f_n(z)t^n$ can be written as
 \begin{equation*}
  P_m(\partial_t,\partial_z)\partial_t^{-m}\left(\sum_{n=0}^\infty u_n(z)t^n\right)+
  Q(\partial_t,\partial_z)\partial_t^{-m}\left(\sum_{n=0}^\infty u_n(z)t^n\right)=\sum_{n=0}^{\infty}f_n(z)t^n.
 \end{equation*}
 
 Since $Q(\partial_t,\partial_z)$ is the rest of the operator $P(\partial_t,\partial_z)$, the lower ordinate of $N(Q)$ is greater than $-m$. Hence, we can write
 $$Q(\partial_t,\partial_z)\partial_t^{-m}\left(\sum_{n=0}^\infty u_n(z)t^n\right)=\sum_{n=1}^\infty\left(\sum_{k=1}^{n}\tilde{Q}(n,k,\partial_z) u_{n-k}(z)\right)t^n$$
for some uniquely defined differential operators $\tilde{Q}(n,k,\partial_z)$, where $n,k\in\N$ and $k\leq n$.
 
 Since also
 \begin{equation*}
 P_m(\partial_t,\partial_z)\partial_t^{-m}\left(\sum_{n=0}^\infty u_n(z)t^n\right)=\sum_{n=0}^\infty 
 \tilde{P}_m(n,\partial_z)u_n(z)t^n,
 \end{equation*}
 where $\tilde{P}_m(n,\partial_z)$ is given by (\ref{eq:operator_p_m_tilde}), we conclude that
 \begin{equation}
 \label{eq2}
  \left\{
   \begin{array}{ll}
    \tilde{P}_m(n,\partial_z)u_n(z)=f_n(z)& \textrm{for}\ n=0\\
    \tilde{P}_m(n,\partial_z)u_n(z)=f_n(z)-\sum_{k=1}^{n}\tilde{Q}(n,k,\partial_z)u_{n-k}(z)& \textrm{for}\ n\geq 1.
   \end{array}
  \right.
 \end{equation}

The fact that the operator $P(\partial_t,\partial_z)\partial_t^{-m}$ is a linear automorphism on $\C[[t,z]]$ means that for every $f=\sum_{n=0}^\infty f_n(z)t^n$ in $\C[[t,z]]$  there exists exactly one $u=\sum_{n=0}^\infty u_n(z) t^n$ in $\C[[t,z]]$ satisfying $P(\partial_t,\partial_z)\partial_t^{-m}u=f$. 

Now we are ready to prove all three implications.
\smallskip\par
((ii), (iii) $\Rightarrow$ (i)) If conditions (ii) and (iii) hold then by Theorem \ref{thm1} the operator $\tilde{P}_m(n,\partial_z)$ is a linear automorphism on $\C[[z]]$ for every $n\in\N_0$. From (\ref{eq2}) it further follows that for every $f=\sum_{n=0}^{\infty}f_n(z)t^n\in\C[[t,z]]$ there exists exactly one $u=\sum_{n=0}^{\infty}u_n(z)t^n\in\C[[t,z]]$ satisfying
\linebreak
$P(\partial_t,\partial_z)\partial_t^{-m}u=f$. Hence, 
$P(\partial_t,\partial_z)\partial_t^{-m}$ is a linear automorphism on $\C[[t,z]]$. 
\smallskip\par
((i), (ii) $\Rightarrow$ (iii)) Assume that the non-resonance condition (iii) does not hold. Then there exists $(n_0,k_0)\in\N_0\times\N_0$ such that $W_{m,0}(n_0,k_0,0)=0$ and $W_{m,0}(n,k,0)\neq0$ for every $n<n_0$ and $k\in\N_0$, and moreover $W_{m,0}(n_0,k,0)\neq0$ for every $k<k_0$. We will show that
$f(t,z)=z^{k_0}t^{n_0}\not\in{\rm im}(P(\partial_t,\partial_z)\partial_t^{-m}, \C[[t,z]])$, and in consequence $P(\partial_t,\partial_z)\partial_t^{-m}$ is not a surjection.
Observe that for any fixed $n<n_0$ the lower ordinate of the Newton polygon $N(\tilde{P}_m(n,\partial_z)$ is equal to zero and the non-resonance condition ($W_{m,0}(n,k,0)\neq 0$ for $k\in\N_0$) for the operator $\tilde{P}_m(n,\partial_z)$ holds. Hence, by Theorem \ref{thm1} the operator $\tilde{P}_m(n,\partial_z)$ given by (\ref{eq:operator_p_m_tilde}) is a linear automorphism on $\C[[z]]$. It means that the only possible solutions $u_n$ ($n<n_0$) of (\ref{eq2}) with $f_n(z)\equiv 0$ for $n<n_0$ are $u_0(z)\equiv 0$, ..., $u_{n_0-1}(z)\equiv 0$. Hence. applying (\ref{eq2}) for $n=n_0$ and $u_0(z)=\dots= u_{n_0-1}(z)=0$, we get $\tilde{P}_m(n_0,\partial_z)u_{n_0}(z)=f_{n_0}(z)=z^{k_0}$. As in the proof of Theorem \ref{thm1} we will show that there exist no $u_{n_0}(z)=\sum_{k=0}^{\infty}u_{n_0,k}z^k$ satisfying the above equation. Indeed, since $W_{m,0}(n_0,k,0)\neq 0$ for every $k<k_0$, as in the proof of Theorem \ref{thm1} we conclude that $u_{n_0,k}=0$ for every $k<k_0$.
Therefore, for $k=k_0$ we get 
$$
W_{m,0}(n_0,k_0,0)u_{n_0,k_0}z^{k_0}=z^{k_0}-\sum_{j=1}^{k_0}\tilde{w}_{k_0-j}u_{n_0,k_0-j}z^{k_0}=z^{k_0},
$$
so $W_{m,0}(n_0,k_0,0)u_{n_0,k_0}=1$, contrary to the condition $W_{m,0}(n_0,k_0,0)=0$.
Thus, $f(t,z)=z^{k_0}t^{n_0}\not\in{\rm im}(P(\partial_t,\partial_z)\partial_t^{-m}, \C[[t,z]])$.
\smallskip\par
((i), (iii) $\Rightarrow$ (ii)) We consider the first equation from the system (\ref{eq2}):
$$
\tilde{P}_m(0,\partial_z)u_0(z)=f_0(z).
$$
Since $P(\partial_t,\partial_z)\partial_t^{-m}$ is an automorphism on $\C[[t,z]]$, we see by (\ref{eq2}) that the operator $\tilde{P}_m(0,\partial_z)$ must be a surjection. On the other hand, since $W_{m,l}(0,k,0)\neq 0$ for every $k\in\N_0$, we conclude that 
$\tilde{P}_m(0,\partial_z)$ is an injection, too. Hence, $\tilde{P}_m(0,\partial_z)$ is an automorphism on $\C[[z]]$. In particular the index $\chi(\tilde{P}_m(0,\partial_z),\C[[z]])$ is equal to zero. On the other hand, by \cite[Corollary 4.2.5]{LodayRichaud} $\chi(\tilde{P}_m(0,\partial_z),\C[[z]])=-l_0$, where $l_0$ is the lower ordinate of the Newton polygon $N(\tilde{P}_m(0,\partial_z))$ defined by (\ref{eq:l}). Thus $l_0=0$ and
$l\leq 0$. Observe also that for sufficiently large $n$, say $n>\tilde{n}_0$, the lower ordinate of the Newton polygon $N(\tilde{P}_m(n,\partial_z))$ is equal to $l$.
Hence, if $l<0$ then $\dim\ker(\tilde{P}_m(n,\partial_z),\C[[z]])>0$ for every $n>\tilde{n}_0$. Moreover, by the non-resonace condition (iii) $\dim\coker(\tilde{P}_m(n,\partial_z),\C[[z]])=0$
(see Remark \ref{re:coker}). In this case the system (\ref{eq2}) has infinitely many solutions, contrary to (i). It means that $l=0$. 

To finish the proof it is sufficient to show that the lower ordinate $l_n$ of the Newton polygon $N(\tilde{P}_m(n,\partial_z))$ is equal to $l$ for every $n\in\N_0$.
Let us suppose the contrary and let $n_1$ be the smallest natural number such that $l_{n_1}\neq l$.
Then by the definitions of $l$ and $l_{n_1}$ we conclude that $l<l_{n_1}$. Hence, 
\begin{equation}
\label{eq:l_n_1}
l_{n_1}>0\quad \textrm{and}\quad {\rm im}(\tilde{P}_m(n_1,\partial_z),\C[[z]])\subseteq z^{l_{n_1}}\C[[z]].
\end{equation}
We will show that in this case $f(t,z)=t^{n_1}\not\in{\rm im}(P(\partial_t,\partial_z)\partial_t^{-m},\C[[t,z]])$. Indeed, since $\tilde{P}_m(n,\partial_z)$ are automorphisms for every $n<n_1$ we conclude that the unique solutions of (\ref{eq2}) with $f_0(z)=\dots=f_{n_1-1}(z)=0$ are $u_0(z)=\dots=u_{n_1-1}(z)=0$. Using (\ref{eq2}) for $n=n_1$
we get the equation $\tilde{P}_m(n_1,\partial_z)u_{n_1}(z)=1$, which has no solutions by the condition (\ref{eq:l_n_1}). Therefore $P(\partial_t,\partial_z)\partial_t^{-m}$ is not a linear automorphism, contrary to (i). Hence, the lower ordinate $l_n$ of the Newton polygon $N(\tilde{P}_m(n,\partial_z))$ is equal to $l$ for every $n\in\N_0$ and $l=0$. It means that (ii) holds.
\end{proof}

\begin{example}\label{ex9}
Let $m\ge0$. We consider the differential operator 
$$P(\partial_t,\partial_z)=p_0(t,z)\partial_{t}^{m}+p_1(t,z)\partial_t^{m+1}+p_2(t,z)\partial_t^{m+1}\partial_z,$$
where 
$$p_0(t,z)=a+zp_{00}(z)+\sum_{n\ge1}p_{0n}(z)t^n\in\C[[t,z]],$$
for some $a>0$ and $p_{0n}\in\C[[z]]$ for all $n\ge0$,
$$p_1(t,z)=(b+zp_{11}(z))t+\sum_{n\ge 2}p_{1n}(z)t^n\in\C[[t,z]],$$
with $b>0$ and $p_{1n}\in\C[[z]]$ for all $n\ge1$, and
$$p_2(t,z)=(c+zp_{21}(z))zt+\sum_{n\ge 2}p_{2n}(z)t^n\in\C[[t,z]],$$
with $c>0$ and $p_{2n}\in\C[[z]]$ for all $n\ge1$. Observe that $\Lambda=\{(m,0),(m+1,0),(m+1,1)\}$ and $\Lambda_{m}=\Lambda$. One has that
$$\tilde{P}_m(n,\partial_z)=a+zp_{00}(z)+(b+zp_{11}(z))n+(c+zp_{21}(z))n\partial_z.$$
If $n=0$, then $\tilde{P}_m(0,\partial_z)=a$ which satisfies that $N(\tilde{P}_m(0,\partial_z))=\corner (0,0)$ whereas $N(\tilde{P}_m(0,\partial_z))=\corner (1,0)$ for $n\ge1$. Therefore, condition (ii) of Theorem~\ref{thm3} is satisfied. In addition to this, $l=0$ and $\Lambda_{m,0}=\Lambda$. We also have $W_{m,0}(n,k,0)=a+bn+cnk\neq0$ for all $(n,k)\in\N_0^2$. Condition (iii) in Theorem~\ref{thm3} holds, and the operator $P(\partial_t,\partial_z)\partial_{t}^{-m}$ is a linear automorphism on $\C[[t,z]]$. 
\end{example}

\begin{remark}
{\rm
The condition ((ii),(iii) imply (i)) from Theorem \ref{thm3}, under which $P(\partial_t,\partial_z)\partial_t^{-m}$ is a~linear automorphism on the space $\C[[t,z]]$, can also be found in \cite[Proposition 1]{LastraTahara} in a more general setting involving nonlinear operators.
}
\end{remark}

\begin{remark}
{\rm 
Theorem \ref{thm3} remains valid when replacing $a_{jr}(t,z)\in\Oo[[t]]$ in (\ref{eq:operator_p}) by $a_{jr}(t,z)\in\C[[t,z]]$ for $(j,r)\in\Lambda$.
}
\end{remark}

\begin{remark}
 {\rm
 If the operator $P(\partial_t,\partial_z)$ is equal to its pricipial part $P_m(\partial_t,\partial_z)$ defined by (\ref{eq:operator_p_m}) then the condition (i) in Theorem \ref{thm4} implies conditions (ii) and (iii).
 
 Indeed, in this case $Q(\partial_t,\partial_z)\equiv 0$ and (\ref{eq2}) is reduced to the system
 \begin{equation}
 \label{eq3}
  \left.
   \begin{array}{ll}
    \tilde{P}_m(n,\partial_z)u_n(z)=f_n(z)& \textrm{for}\ n\geq 0.
   \end{array}
  \right.
 \end{equation}
 
 The condition (i) means that for every $\sum_{n=0}^{\infty}f_n(z)t^n\in\C[[t,z]]$ the exists exactly one $\sum_{n=0}^{\infty}u_n(z)t^n\in\C[[t,z]]$ such that (\ref{eq3}) holds. Hence, for every $f_n(z)\in\C[[z]]$ there exists exactly one $u_n(z)\in\C[[z]]$ such that
 (\ref{eq3}) is satisified for $n\geq 0$. Thus the operator $\tilde{P}_m(n,\partial_z)$ is a~linear automorphism on $\C[[z]]$ for every $n\in\N_0$. Applying Theorem \ref{thm1} to the operator $\tilde{P}_m(n,\partial_z)$ for every $n\in\N_0$ we conclude that the conditions (ii) and (iii) in Theorem \ref{thm3} are satisfied.
 }
\end{remark}

\begin{remark}
 {\rm
 It is still an open problem whether in Theorem \ref{thm3} the condition (i) implies the conditions (ii) and (iii) for any operator $P(\partial_t,\partial_z)$ defined by (\ref{eq:operator_p}).
 }
\end{remark}

Theorem~\ref{thm3} can be extended to the more general framework of Section~\ref{sec31} without difficulty. Indeed, assume that $\mathfrak{m}_1=(m_{1,n})_{n\ge0}$ and $\mathfrak{m}_2=(m_{2,n})_{n\ge0}$ are two sequences of positive real numbers normalized by $m_{1,0}=m_{2,0}=1$. 

We define the operator $\partial_{\mathfrak{m}_1,t}^{-1}$ as the formal inverse operator of $\partial_{\mathfrak{m}_1,t}$, whose images do not have constant term. $\partial_{\mathfrak{m}_1,t}^{-m}$ for $m\in\N_0$ is defined recursively in a natural way.

Let us consider a moment integro-differential operator of the form 
$$P(\partial_{\mathfrak{m}_1,t},\partial_{\mathfrak{m}_{2},z})\partial_{\mathfrak{m}_1,t}^{-m}:\C[[t,z]]\to\C[[t,z]],$$ with $P$ as in (\ref{eq:operator_p}) substituting classical derivatives by their corresponding moment differential operators. All the constructions can be adapted to this framework in a straightforward manner. Indeed, such construction provides
$$\tilde{P}_m(n,\partial_{\mathfrak{m}_2,z})=\sum_{(j,r)\in\Lambda_m}\tilde{a}_{j,r,j-r}(z)z^{\alpha_{j,r}}\frac{m_{1,n}}{m_{1,n-(j-m)}}\partial_{\mathfrak{m}_2,z}^r$$
as the generalization of (\ref{eq:operator_p_m_tilde}), and 
 \begin{equation*}
 W_{m,l,\mathfrak{m}_1,\mathfrak{m}_2}(\lambda_1,\lambda_2,z)=\sum_{\genfrac{}{}{0pt}{2}{(j,r)\in\Lambda_{m,l}}{j\le \lambda_1+m,\, r\le \lambda_2}}\tilde{a}_{j,r,j-m}(z)\frac{m_{1,\lambda_1}m_{2,\lambda_2}}{m_{1,\lambda_1-(j-m)}m_{2,\lambda_2-r}}
 \end{equation*}
generalizing (\ref{eq:W_ml}). Theorem~\ref{thm3} reads as follows in this context.

\begin{corol}
Any two of the following three conditions entail a third one:
  \begin{enumerate}
 \item[(i)] The operator $P(\partial_{\mathfrak{m}_1,t},\partial_{\mathfrak{m}_2,z})\partial_{\mathfrak{m}_1,t}^{-m}$ is a linear automorphism on $\C[[t,z]]$.
  \item[(ii)] The lower ordinate of the Newton polygon $N(\tilde{P}_m(n,\partial_{\mathfrak{m}_2,z}))$ is equal to zero for every $n\in\N_0$.
  \item[(iii)] The non-resonance condition 
  $W_{m,l,\mathfrak{m}_1,\mathfrak{m}_2}(n,k,0)\neq0$  holds for every $n,k\in\N_0$.
 \end{enumerate}
\end{corol}

At this point, the results in Section~\ref{sec32} can also be applied to achieve a result under Gevrey settings in the several variable case. Here we assume that the operator $P(\partial_t,\partial_z)$ defined by (\ref{eq:operator_p}) has coefficients $a_{jr}(t,z)\in\Oo[[t]]$.

\begin{theo}
 \label{thm4}
 Let $s\geq 0$. Assume that the following conditions hold:
 \begin{enumerate}
 \item[(a)] The lower ordinate $l$ of the Newton polygon $N(\tilde{P}_m(n,\partial_z))$ is equal to zero for every $n\in\N_0$.
  \item[(b)] The first positive slope of the Newton polygon $N(\tilde{P}_m(n,\partial_z))$ is greater or equal to $1/s$ for every $n\in\N_0$.
  \item[(c)] The non-resonance condition $W_{m,l}(n,k,0)\ne 0$ holds for every $n,k\in\N_0$.
 \end{enumerate}
 Then the operator $P(\partial_t,\partial_z)\partial_t^{-m}$ is a linear automorphism on $\C[[t,z]]$, which extends to an automorphism on $\C[[z]]_s[[t]]$.
 
 Conversely, if the operator $P(\partial_t,\partial_z)\partial_t^{-m}$ is a linear automorphism on $\C[[t,z]]$, which extends to an automorphism on $\C[[z]]_s[[t]]$, and one of the conditions (a) or (c) is true, then all conditions (a), (b) and (c) are satisified. 
\end{theo}
\begin{proof}

($\Rightarrow$) Assume that statements (a), (b) and (c) hold. 
By Theorem \ref{thm3} the operator $\tilde{P}_m(n,\partial_z)$ is a linear automorphism on $\C[[z]]$ for every $n\in\N_0$. Theorem \ref{thm2} guarantees that the operator $\tilde{P}_m(n,\partial_z)$, as a linear automorphism on $\C[[z]]$, extends to a linear automorphism on $\C[[z]]_s$ for every $n\in\N_0$.

This also entails that $P(\partial_t,\partial_z)\partial_t^{-m}$ is a linear automorphism on $\C[[t,z]]$ which extends to the linear automorphism on $\C[[z]]_s[[t]]$. 

($\Leftarrow$) By Theorem \ref{thm3}, for every $n\in\N_0$ the operator $\tilde{P}_m(n,\partial_z)$ is a linear automorphism on $\C[[z]]$, which extends by the assumption to the automorphism on $\C[[z]]_s$. 

Hence, by Theorem \ref{thm2} we conclude that the lower ordinate $l$ of the Newton polygon $N(\tilde{P}_m(n,\partial_z))$ is equal to zero and its first positive slope is greater or equal to $1/s$. We also observe that $W_{m,0}(n,k,0)\ne 0$ for every $k\in\N_0$. The previous statement holds for every $n\in\N_0$, so conditions (a), (b) and (c) in Theorem \ref{thm4} are satisified.
\end{proof}

As a direct consequence of Theorem \ref{thm4} with $s=0$ we get
\begin{corol}
\label{cor:3}
 Assume that the following conditions hold:
 \begin{enumerate}
  \item[(a)] $N(\tilde{P}_m(n,\partial_z))=\corner(p_n,0)$, where $p_n$ denotes the order of the operator $\tilde{P}_m(n,\partial_z)$ for every $n\in\N_0$,
  \item[(b)] The non-resonance condition $W_{m,l}(n,k,0)\ne 0$ holds for every $n,k\in\N_0$.
 \end{enumerate}
 Then the operator $P(\partial_t,\partial_z)\partial_t^{-m}$ is a linear automorphism on $\C[[t,z]]$ which extends to an automorphism on $\C[[z]]_0[[t]]$.
 
 Conversely, if the operator $P(\partial_t,\partial_z)\partial_t^{-m}$ is a linear automorphism on $\C[[t,z]]$, which extends to an automorphism on $\C[[z]]_0[[t]]$ and one of the conditons (a) or (b) holds, then the remaining condition also is satisified.
\end{corol}

\begin{remark}
{\rm By the corollary above we get the sufficient condition under which the operator $P(\partial_t,\partial_z)\partial^{-m}_t$ is a linear automorphism on the space $\C[[z]]_0[[t]]$ of formal power series in $t$ with coefficients being the germs of holomorphic functions in $z$ at zero. Observe that $u(t,z)=\sum_{n=0}^{\infty}u_n(z)t^n\in\C[[z]]_0[[t]]$ if and only if there exists a sequence of positive real numbers $(r_n)_{n\geq0}$ such that $u_n\in\Oo(D(r_n))$ for $n\in\N_0$, where $D(r_n)$ denotes the complex disc in $\C$ with center at zero and radius $r_n$. In particular, taking the constant series $(r_n)_{n\geq0}=(r)_{n\geq0}$ for some $r>0$ we get the series $u(t,z)$ belonging to the space $\Oo[[t]]$ being the subspace of $\C[[z]]_0[[t]]$.}

{\rm In the subsequent paper we are planning to study the relation between the sequences $(r_n)_{n\geq 0}$ describing the solution $u(t,z)$, and the operator $P(\partial_t,\partial_z)\partial_t^{-m}$.
}

{\rm It is worth emphasizing that recently there has been a growing interest in such spaces of formal power series with coefficients being holomorphic on different complex neighborhood of the origin, see \cite{carrillolastra} in the formal settings, and \cite{tahara3} and \cite{lama4} regarding summability results. 
}
\end{remark}

\begin{example}
Observe that the differential operator in Example~\ref{ex9} can be modified by adding any number of terms of the form $a_{jr}(t,z)\partial_t^j\partial_z^r$, with $j-\hbox{ord}_t(a_{jr})<m$. In this situation, $\Lambda_m$ coincides with that in Example~\ref{ex9}, so $m$ and $\tilde{P}_m(n,\partial_z)$ are the same for both differential operators, and the conditions of Theorem~\ref{thm3} are also satisfied.
\end{example}

As a conclusion of the present work, we link our results to previous works which can be found in the literature. 

Using Corollary \ref{cor:3} we get the known sufficient conditions under which the operator $P(\partial_t,\partial_z)\partial_t^{-m}$ is a linear automorphism on the space $\C[[t,z]]$ or $\C[[z]]_0[[t]]$. A first example is the following weaker version of the results for the space $\Oo[[t]]$, which one can find in  \cite[Proposition 1.1]{BalserLodayRichaud}, \cite[Theorem 1]{Remy1} or \cite[Theorem 1]{Remy2}.

\begin{corol}
Let $(m,-m)$ be a vertex of $N(P)$ (i.e., $\Lambda_m=\{(m,0)\}$, $\ord_t(a_{m,0})=0$ and $P_m(\partial_t,\partial_z)=a_{m,0,0}(z)\partial_t^m$) and assume the non-resonance condition $a_{m,0,0}(0)\neq 0$ holds. Then, the operator $P(\partial_t,\partial_z)\partial_t^{-m}$ is a linear automorphism on $\C[[z]]_0[[t]]$.
\end{corol}

More generally we have the following result. This type of results for the space $\Oo[[t]]$ are also studied in \cite[Theorem 1]{BaouendiGoulaouic}.

\begin{corol}
Let $\Lambda_m\subset \N_0\times\{0\}$ (i.e., $P_m(\partial_t,\partial_z)=\sum_{(j,0)\in\Lambda_m}a_{j,0,j-m}(z)t^{j-m}\partial_t^j$) such that its associated characteristic polynomial
 $$W(\lambda,z):=\sum_{(j,0)\in\Lambda_m}a_{j,0,j-m}(z)\lambda(\lambda-1)\cdots(\lambda-(j-m)+1)$$
 satisfies the non-resonance condition $W(n,0)\neq 0$ for every $n\in\N_0$. Then, the Fuchsian operator $P(\partial_t,\partial_z)\partial_t^{-m}$ is a linear automorphism on $\C[[z]]_0[[t]]$.
\end{corol}

\end{document}